\documentclass[12pt,reqno]{amsart}
\usepackage{amsmath, amsfonts, amssymb, amsthm, hyperref}
\usepackage{bm}
\allowdisplaybreaks[4]
\def\l{\left}
\def\r{\right}
\def\bg{\bigg}
\def\({\bg(}
\def\){\bg)}
\def\t{\text}
\def\f{\frac}

\def\sgn{{\rm sign}}

\def\ls{\leqslant}
\def\gs{\geqslant}

\def\sm{\setminus}
\def\bi{\binom}

\def\Fix{{\rm Fix}}

\def\Proof{\noindent{\it Proof}}

\def\N{\mathbb N}

\def\1{{\bf 1}}

\def\pmod #1{\ ({\rm{mod}}\ #1)}

\def\per{{\rm per}}

\def\<{\langle}
\def\>{\rangle}
\theoremstyle{plain}
\newtheorem{theorem}{Theorem}[section]
\newtheorem{lemma}{Lemma}

\newtheorem{conjecture}{Conjecture}

\theoremstyle{definition}

\theoremstyle{remark}
\newtheorem{remark}{Remark}

\makeatletter
\@namedef{subjclassname@2020}{%
  \textup{2020} Mathematics Subject Classification}
\makeatother
 \vspace{4mm}

\begin{document}
\hbox{Preprint, {\tt arXiv:2208.12167}}
\medskip

\title[A novel permanent identity with applications]
{A novel permanent identity with applications}
\author{Yue-Feng She}
\address {(Yue-Feng She) Department of Mathematics, Nanjing
University, Nanjing 210093, People's Republic of China}
\email{she.math@smail.nju.edu.cn}

\author{Zhi-Wei Sun}
\address {(Zhi-Wei Sun, corresponding author) Department of Mathematics, Nanjing
University, Nanjing 210093, People's Republic of China}
\email{zwsun@nju.edu.cn}

\author{Wei Xia}
\address {(Wei Xia) Department of Mathematics, Nanjing
	University, Nanjing 210093, People's Republic of China}
\email{wxia@smail.nju.edu.cn}

\keywords{Permanent, permutation, combinatorial identity, rational function.
\newline \indent Supported by the Natural Science Foundation of China (grant no. 11971222).}
\begin{abstract}
Let $n$ be a positive integer, and define the rational function $S(x_1,\ldots,x_{2n})$ as the permanent of the matrix $[x_{j,k}]_{1\ls j,k\ls 2n}$, where
$$x_{j,k}=\begin{cases}(x_j+x_k)/(x_j-x_k)&\text{if}\ j\not=k,\\1&\text{if}\ j=k.\end{cases}$$
We give an explicit formula for $S(x_1,\ldots,x_{2n})$ which has the following consequence: If one of the variables $x_1,\ldots,x_{2n}$ takes zero, then
$S(x_1,\ldots,x_{2n})$ vanishes, i.e.,
$$\sum_{\tau\in S_{2n}}\prod_{j=1\atop \tau(j)\not=j}^{2n}\f{x_j+x_{\tau(j)}}{x_j-x_{\tau(j)}}=0,$$
where we view an empty product $\prod_{i\in\emptyset}a_i$ as $1$. As an application, we show that
if $\zeta$ is a primitive $2n$-th root of unity then
$$\sum_{\tau\in S_{2n}}\prod_{j=1\atop \tau(j)\not=j}^{2n}\f{1+\zeta^{j-\tau(j)}}{1-\zeta^{j-\tau(j)}}=((2n-1)!!)^2$$
as conjectured by Z.-W. Sun.
\end{abstract}
\maketitle

\section{Introduction}
\setcounter{lemma}{0}
\setcounter{theorem}{0}
\setcounter{equation}{0}
\setcounter{conjecture}{0}
\setcounter{remark}{0}
\setcounter{corollary}{0}

A {\it permanent} of an $n\times n$ matrix $A=[a_{j,k}]_{1\ls j,k\ls n}$ over a field, is defined as
$$\per(A)=\per[a_{j,k}]_{1\ls j,k\ls n}=\sum_{\tau\in S_n}\prod_{j=1}^na_{j,\tau(j)},$$
where $S_n$ is the symmetric group consisting of all the permutations of $\{1,\ldots,n\}$.
Z.-W. Sun \cite{S1} investigated arithmetic properties of some permanents.
Motivated by this, S. Fu, Z. Lin and Sun \cite{FLS} proved several permanent identities related to some combinatorial sequences like the Bernoulli numbers. Sun \cite[\S11.3]{ConjBook} contains several  conjectures involving permanents.

For any positive integer $n$, we define the rational function $S(x_1,\ldots,x_{2n})$ as
$\per[x_{j,k}]_{1\ls j,k\ls 2n}$, where
\begin{equation}\label{xjk}x_{j,k}=\begin{cases}(x_j+x_k)/(x_j-x_k)&\text{if}\ j\not=k,\\1&\text{if}\ j=k.\end{cases}
\end{equation}
In other words,
\begin{equation}\label{Def-S} S(x_1,\ldots,x_{2n})=\sum_{\tau\in S_{2n}}\prod_{j=1\atop \tau(j)\not=j}^{2n}\f{x_j+x_{\tau(j)}}{x_j-x_{\tau(j)}},
\end{equation}
where we view an empty product $\prod_{i\in\emptyset}a_i$ as $1$.
For example,
\begin{equation}\label{S(x12)} S(x_1,x_2)=\sum_{\tau\in S_2}\prod_{j=1\atop \tau(j)\not=j}^2\f{x_j+x_{\tau(j)}}{x_j-x_{\tau(j)}}
=1+\f{x_1+x_2}{x_1-x_2}\times\f{x_2+x_1}{x_2-x_1}=1-\l(\f{x_1+x_2}{x_1-x_2}\r)^2,
\end{equation}
which vanishes if $x_1$ or $x_2$ is zero.

The second author posed the following conjecture
in \cite[A356041]{OEIS}.

\begin{conjecture}\label{Conj} Let $n$ be any positive integer.
If one of the variables $x_1,\ldots,x_{2n}$ takes zero, then
\begin{equation*}  S(x_1,\ldots,x_{2n})=0,\ \t{i.e.},\
\sum_{\tau\in S_{2n}}\prod_{j=1\atop \tau(j)\not=j}^{2n}\f{x_j+x_{\tau(j)}}{x_j-x_{\tau(j)}}=0.
\end{equation*}
\end{conjecture}

\begin{remark} Conjecture 1.1 is somewhat surprising, and it was motivated by the following open conjecture of the second author \cite{S2}: For any prime $p$, we have
$$\sum_{\tau\in S_{p-1}}\prod_{j=1\atop \tau(j)\not=j}^{p-1}\frac{j+\tau(j)}{j-\tau(j)}\equiv((p-2)!!)^2\pmod{p^2}.$$
\end{remark}

In this paper, we will confirm Conjecture \ref{Conj} and establish the following more general result.

\begin{theorem}\label{Th1.1} For any positive integer $n$ and the rational function $S(x_1,\ldots,x_{2n})$ defined by \eqref{Def-S}, we have
	\begin{equation} \label{S} S(x_1,\ldots, x_{2n})=(-4)^nx_1x_2\cdots x_{2n}\sum_{\bigcup_{k=1}^n\{i_k,j_k\}=\{1,\ldots,2n\}} \prod_{k=1}^n\frac{1}{(x_{i_k}-x_{j_k})^2},
	\end{equation}
where the sum runs over all partitions of $\{1,\ldots,2n\}$ into $n$ pairs.
\end{theorem}

In contrast with Theorem 1.1, it is interesting to investigate the sequence $r_n=\per[m_{j,k}]_{0\ls j,k\ls 2n}$ $(n=1,2,3,\ldots)$, where
$$m_{j,k}=\begin{cases}(j+k)/(j-k)&\t{if}\ j\not=k,\\1&\t{if}\ j=k.\end{cases}$$
The initial four terms of this sequence (cf. \cite[A356041]{OEIS}) are
$$r_1=-10,\ r_2=\f{5870}9,\ r_3=-\f{436619903}{4050},\ r_4=\f{204409938157631}{6125000}.$$
We are unable to evaluate $r_n$ for a general positive integer $n$.

Recall that a cycle (or cyclic permutation) in $S_{2n}$ of length $l$ has the form $(a_1\ \ldots\ a_l)$,
 where $a_1,\ldots,a_l$ are distinct elements of $\{1,\ldots,2n\}$. A cycle of length $l$ is simply called an $l$-cycle.

For each $k\in\{1,2,\dots,2n\}$, we define
$$C(k)=\{\tau\in S_{2n}:\ \tau\ \t{is a}\ k\t{-cycle}\}.$$
In contrast with the equality \eqref{Def-S}, we also set
\begin{equation}\label{Def-s}s(x_1,\ldots,x_{2n})=\sum_{\tau\in C(2n)}\prod_{j=1\atop \tau(j)\not=j}^{2n}\f{x_j+x_{\tau(j)}}{x_j-x_{\tau(j)}}.\end{equation}
For example,
$$s(x_1,x_2)=f((1\,2))=\f{x_1+x_2}{x_1-x_2}\times\f{x_2+x_1}{x_2-x_1}
=-\left(\frac{x_1+x_2}{x_1-x_2}\right)^2.$$

We find that $s(x_1,\ldots,x_{2n})$
is related to the tangent numbers $T_n\ (n\in\N)$ defined by
$$\tan x=\sum_{n=1}^\infty T_n\f{x^{2n-1}}{(2n-1)!}\ \l(|x|<\f{\pi}2\r).$$
 For formulas and combinatorial interpretations of the tangent numbers, one may consult
Sloane \cite[A000182]{Sloane}. For example, it is known that 
\begin{equation} (-1)^nT_n=2^{2n}(1-2^{2n})\f{B_{2n}}{2n}\ \quad\t{for all}\ n=1,2,3,\ldots,
\end{equation}
where the Bernoulli numbers $B_0,B_1,B_2,\ldots$ are given by
$$\f x{e^x-1}=\sum_{n=0}^\infty B_n\f{x^n}{n!}\ \l(|x|<2\pi\r).$$

Now we state our second theorem.

\begin{theorem}\label{Th1.2}
	If $n>1$ is an integer, then
	\begin{equation}\label{s-T}
		s(x_1,\ldots,x_{2n})=(-1)^nT_n.
	\end{equation}
\end{theorem}
\begin{remark} Theorem \ref{Th1.2} plays an important role in our proof of Theorem \ref{Th1.1}.
\end{remark}

Via applying Theorem \ref{Th1.1}, we prove the following result part (ii) of which was first conjectured by the second author \cite{Sun}.

\begin{theorem}\label{Th1.3} {\rm (i)} Let $n$ be any positive integer. Then
\begin{equation} \label{per-det identity}
	\per[x_{j,k}-1]_{1\ls j,k\ls 2n}=(-1)^n\det[x_{j,k}-1]_{1\ls j,k\ls 2n}=S(x_1,\ldots,x_{2n}),
\end{equation}
where $x_{j,k}$ is given by \eqref{xjk}.

{\rm (ii)} Let $n>1$ be an integer, and let $\zeta$ be a primitive $n$th root of unity. Set
	\begin{equation*}
		c_{j,k}=\begin{cases} \f{1+\zeta^{j-k}}{1-\zeta^{j-k}} &\t{if $j\not=k$,}\\1 &\t{if $j=k$.}\end{cases}			
	\end{equation*}
If $n$ is even, then
\begin{equation}\label{c-even}\per[c_{j,k}]_{1\ls j,k\ls n}=((n-1)!!)^2.
\end{equation}
When $n$ is odd, we have
\begin{equation}\label{c-odd}\per[c_{j,k}]_{1\ls j,k\ls n-1}=\f{((n-1)!!)^2}n.
\end{equation}
\end{theorem}
\begin{remark} Let $n>1$ be an integer. By H. Wang and Sun \cite{WangSun}, we have
$$\det[c_{j,k}]_{1\ls j,k\ls n}=\prod_{s=1}^{n-1}(n+1-2s)=\begin{cases}(-1)^{n/2-1}((n-1)!!)^2/(n-1)
&\t{if}\ 2\mid n,\\0&\t{if}\ 2\nmid n.\end{cases}$$
When $n$ is odd, Wang and Sun \cite{WangSun} also proved that
$$\det[c_{j,k}]_{1\ls j,k\ls n-1}=(-1)^{(n+1)/2}\f{((n-1)!!)^2}{n(n-1)}.$$
\end{remark}

In the next section we provide three lemmas. Section 3 is devoted to our proof of Theorem \ref{Th1.2}. We will prove Theorems \ref{Th1.1} and \ref{Th1.3} in Sections 4 and 5 respectively.

\section{Three Lemmas}
\setcounter{lemma}{0}
\setcounter{theorem}{0}
\setcounter{equation}{0}
\setcounter{conjecture}{0}
\setcounter{remark}{0}
\setcounter{corollary}{0}

As usual, for a set $X$ we let $\#X$ denote the cardinality of $X$.
If $X$ is a nonempty set, then we let $S(X)$
denote the symmetric group consisting of all permutations of $X$.
If $X=\{1,\ldots,n\}$, then $S(X)$ is usually denoted by $S_n$.

Let $n$ be a positive integer. For any $\tau\in S_{2n}$, we set
$$\Fix(\tau)=\{1\ls j\ls 2n:\ \tau(j)=j\}
\ \ \t{and}\ \ D(\tau)=\{1,\ldots,2n\}\sm\Fix(\tau).$$
If $D(\tau)=\{1,\ldots,2n\}$ (i.e., $\Fix(\tau)=\emptyset$),
then $\tau$ is called a {\it derangement} of $1,\ldots,2n$.

In this section,
we let $x_1, x_2, \dots,x_{2n}$ be variables. For any $\tau\in S_{2n}$, we set
$$f(\tau)=\prod_{j\in D(\tau)}\f{x_j+x_{\tau(j)}}{x_j-x_{\tau(j)}}.$$
Note that $f(\tau)=1$ if $D(\tau)=\emptyset$ (i.e., $\tau$ is the identity of the group $S_{2n}$).

 \begin{lemma}\label{Lem2.1} Let $\sigma,\tau\in S_{2n}$ with $D(\sigma)\cap D(\tau)=\emptyset$.
 Then
 \begin{equation}\label{sigma-tau} \sigma\tau=\tau\sigma,\ D(\sigma\tau)=D(\sigma)\cup D(\tau)\ \t{and}\ f(\sigma\tau)=f(\sigma)f(\tau).
 \end{equation}
 \end{lemma}
 \Proof. (i) Let $j\in\{1,\ldots,2n\}$. We first prove that $\sigma\tau(j)=\tau\sigma(j)$.

 If $j\not\in D(\sigma)\cup D(\tau)$, then $\sigma(j)=j$ and $\tau(j)=j$, thus
 $$\sigma\tau(j)=\sigma(\tau(j))=\sigma(j)=j=\tau(j)=\tau(\sigma(j)).$$

 Now suppose that $j\in D(\sigma)\cup D(\tau)$. Without loss of generality, we assume
 $j\in D(\sigma)$.
 Note that $j\not\in D(\tau)$ and thus $\tau(j)=j$.
 As $\sigma(j)\not=j$, we also have $\sigma(\sigma(j))\not=\sigma(j)$,
 hence $\sigma(j)\in D(\sigma)$ and thus $\sigma(j)\not\in D(\tau)$,
 so
 $$\sigma\tau(j)=\sigma(j)=\tau\sigma(j).$$

 (ii) Next we show that $D(\sigma\tau)=D(\sigma)\cup D(\tau)$, i.e., for any $1\ls j\ls 2n$ we have
 \begin{equation}\label{D(sigma-tau)}\sigma\tau(j)\not=j\iff \sigma(j)\not=j\ \t{or}\ \tau(j)\not=j.
 \end{equation}
 The ``$\Rightarrow$" direction is obvious. If $\sigma(j)\not=j$ (i.e, $j\in D(\sigma)$),
 then $j\not\in D(\tau)$ and hence $\sigma\tau(j)=\sigma(j)\not=j$.
 If $\tau(j)\not=j$, then $j\not\in D(\sigma)$ and hence $\sigma\tau(j)=\tau\sigma(j)=\tau(j)\not=j$.
 Therefore \eqref{D(sigma-tau)} holds.

 (iii) Finally we prove that $f(\sigma\tau)=f(\sigma)f(\tau)$.
 Note that $\tau(j)=j$ if $j\in D(\sigma)$, and $\sigma(j)=j$ if $j\in D(\tau)$.
 In view of (i) and (ii), we have
 \begin{align*}f(\sigma\tau)&=\prod_{j\in D(\sigma)\cup D(\tau)}\f{x_j+x_{\sigma\tau(j)}}{x_j-x_{\sigma\tau(j)}}
 \\&=\prod_{j\in D(\sigma)}\f{x_j+x_{\sigma\tau(j)}}{x_j-x_{\sigma\tau(j)}}
 \times\prod_{j\in D(\tau)}\f{x_j+x_{\tau\sigma(j)}}{x_j-x_{\tau\sigma(j)}}
 \\&=\prod_{j\in D(\sigma)}\f{x_j+x_{\sigma(j)}}{x_j-x_{\sigma(j)}}
 \times\prod_{j\in D(\tau)}\f{x_j+x_{\tau(j)}}{x_j-x_{\tau(j)}}
 \\&=f(\sigma)f(\tau).
 \end{align*}

 Combining the above, we have completed the proof of \eqref{sigma-tau}. \qed

\begin{lemma}\label{Lem2.2} For any $\tau\in S_{2n}$, we have
\begin{equation}\label{Inv}f(\tau^{-1})=(-1)^{\#D(\tau)}f(\tau).
\end{equation}
\end{lemma}
\Proof. Observe that
\begin{align*}f(\tau^{-1})&=\prod^{2n}_{j=1\atop \tau^{-1}(j)\not=j}\f{x_j+x_{\tau^{-1}(j)}}{x_j-x_{\tau^{-1}(j)}}
=\prod^{2n}_{i=1\atop i\not=\tau(i)}\f{x_{\tau(i)}+x_i}{x_{\tau(i)}-x_i}
\\&=(-1)^{|D(\tau)|}\prod^{2n}_{i=1\atop\tau(i)\not=i}\f{x_i+x_{\tau(i)}}{x_i-x_{\tau(i)}}
=(-1)^{|D(\tau)|}f(\tau).
\end{align*}
This proves \eqref{Inv}. \qed

Recall that $s(x_1,x_2,\dots,x_{2n})$ is defined by \eqref{Def-s}.

\begin{lemma}\label{Lem2.3} Let $n>1$ be an integer. Then
$s(x_1,x_2,\dots,x_{2n})$ is a constant only depending on $n$.
\end{lemma}
\Proof. Suppose $\{m_1,m_2,\dots,m_{2n}\}=\{1,2,\dots,2n\}$ and set $m_{2n+1}=m_2$. Define
\begin{align*}
	\sigma_1=\ &(m_1\ m_2\ m_3\ \cdots\ m_{2n-1}\ m_{2n}),\\
	\sigma_2=\ &(m_1\ m_3\ m_4\ \cdots\ m_{2n}\ m_2),\\
	&\vdots\\
	\sigma_{2n-1}=\ &(m_1\ m_{2n}\ m_2\ \cdots\ m_{2n-1}).\\
\end{align*}
Note that
\begin{align*} f(\sigma_1)&=\f{x_{m_1}+x_{m_2}}{x_{m_1}-x_{m_2}}\times\prod_{i=2}^{2n-1}\f{x_{m_{i}}+x_{m_{i+1}}}{x_{m_{i}}
-x_{m_{i+1}}}
\times\f{x_{m_{2n}}+x_{m_1}}{x_{m_{2n}}-x_{m_1}}
\\&=X\times\f{x_{m_{2n}}-x_{m_{2n+1}}}{x_{m_{2n}}+x_{m_{2n+1}}}\times\f{x_{m_1}+x_{m_2}}{x_{m_1}-x_{m_2}}
\times\f{x_{m_{2n}}+x_{m_1}}{x_{m_{2n}}-x_{m_1}}
\\&=X\times\f{x_{m_{2n}}-x_{m_{2n+1}}}{x_{m_{2n}}+x_{m_{2n+1}}}\times\f{x_{m_1}+x_{m_{2n+1}}}
{x_{m_1}-x_{m_{2n+1}}}
\times\f{x_{m_{2n}}+x_{m_1}}{x_{m_{2n}}-x_{m_1}},
\end{align*}
where
$$X=\prod_{j=2}^{2n}\f{x_{m_j}+x_{m_{j+1}}}{x_{m_j}-x_{m_{j+1}}}.$$
For each $i=2,\ldots,2n-1$, as $$\sigma_i=(m_1\ m_{i+1}\ \cdots\ m_{2n}\ m_2\ \cdots m_i),$$ we have
\begin{align*}f(\sigma_i)=\ &\f{x_{m_1}+x_{m_{i+1}}}{x_{m_1}-x_{m_{i+1}}}
\times\f{x_{m_{2n}}+x_{m_2}}{x_{m_{2n}}-x_{m_2}}\times\f{x_{m_i}+x_{m_1}}{x_{m_i}-x_{m_1}}
\\\ &\times\prod_{1<j<i\ \t{or}\ i<j<2n}\f{x_{m_j}+x_{m_{j+1}}}{x_{m_j}-x_{m_{j+1}}}
\\=\ &X\times\f{x_{m_i}-x_{m_{i+1}}}{x_{m_i}+x_{m_{i+1}}}\times\f{x_{m_1}+x_{m_{i+1}}}{x_{m_1}-x_{m_{i+1}}}
\times\f{x_{m_i}+x_{m_1}}{x_{m_i}-x_{m_1}},
\end{align*}
Therefore
$$\sum_{i=1}^{2n-1}f(\sigma_i)=X\sum_{i=2}^{2n}\f{x_{m_i}-x_{m_{i+1}}}{x_{m_i}+x_{m_{i+1}}}
\times\f{x_{m_1}+x_{m_{i+1}}}{x_{m_1}-x_{m_{i+1}}}
\times\f{x_{m_i}+x_{m_1}}{x_{m_i}-x_{m_1}}.$$
It is easy to verify the identity
\begin{equation}\f{y-z}{y+z}\times\f{x+z}{x-z}\times\f{y+x}{y-x}
=\f{z-y}{z+y}+2x\l(\f1{x-z}-\f1{x-y}\r).\end{equation}
So, by the above, we have
\begin{align*}\sum_{i=1}^{2n-1}f(\sigma_i)&=X\sum_{i=2}^{2n}\l(\f{x_{m_{i+1}}-x_{m_{i}}}{x_{m_{i+1}}+x_{m_i}}
+2x_{m_1}\l(\f1{x_{m_1}-x_{m_{i+1}}}-\f1{x_{m_1}-x_{m_i}}\r)\r)
\\&=X\sum_{i=2}^{2n}\f{x_{m_{i+1}}-x_{m_{i}}}{x_{m_{i+1}}+x_{m_i}}
+2x_{m_1}X\l(\f1{x_{m_1}-x_{m_{2n+1}}}-\f1{x_{m_1}-x_{m_2}}\r)
\\&=X\sum_{i=2}^{2n}\f{x_{m_{i+1}}-x_{m_{i}}}{x_{m_{i+1}}+x_{m_i}},
\end{align*}
which is independent of $x_{m_1}$.

By the above, for each $1\ls m_1\ls 2n$,  $s(x_1,\ldots,x_{2n})$
is independent of $x_{m_1}$. So $s(x_1,\ldots,x_{2n})$ is a constant only depending on $n$.
This concludes the proof. \qed

\section{Proof of Theorem 1.2}
\setcounter{lemma}{0}
\setcounter{theorem}{0}
\setcounter{equation}{0}
\setcounter{conjecture}{0}
\setcounter{remark}{0}
\setcounter{corollary}{0}

For each $n=1,2,3,\ldots$, we set
\begin{equation}A_n=[a_{i,j}]_{1\ls i,j\ls 2n},\ \ \t{where}\
	a_{i,j}=\begin{cases}
		1 &\t{if $i\gs j,$}\\
		-1  &\t{if $i<j$.}	\end{cases}
\end{equation}

\begin{lemma}\label{Lem2.3}
	For any positive integer $n$, we have $\per(A_n)=0.$
\end{lemma}
\Proof. Observe that
\begin{align*}\per(A_n)&=\sum_{\sigma\in S_{2n}}\prod_{i=1}^{2n}a_{i,\sigma(i)}
\\&=\sum_{t=1}^{2n}\sum_{\sigma\in S_{2n}\atop\sigma(t)=1}a_{t,\sigma(t)}\prod_{i=1\atop i\not=t}^{2n}
a_{i,\sigma(i)}
\\&=\sum_{t=1}^{2n}a_{t1}\sum_{\sigma\in S_{2n}\atop \sigma(t)=1}\prod_{1\ls i<t}a_{i,\sigma(i)}
\times\prod_{t<i\ls 2n}a_{i,\sigma(i)}
\\&=\sum_{t=1}^{2n}a_{t1}\,\per[a_{i,j}^{(t)}]_{1\ls i,j\ls 2n-1}=\sum_{t=1}^{2n}\per[a_{i,j}^{(t)}]_{1\ls i,j\ls 2n-1},
\end{align*}
where
$$a_{i,j}^{(t)}=\begin{cases}a_{i,j+1}&\t{if}\ 1\ls i<t,
\\a_{i+1,j+1}&\t{if}\ t\ls i\ls 2n.\end{cases}$$
It is easy to see that
$$a_{i,j}^{(t)}=\begin{cases}1&\t{if}\ i>j\ \t{or}\ i=j\gs t,
\\-1&\t{if}\ i<j\ \t{or}\ i=j<t.\end{cases}$$
Therefore, for the matrix $A_n^{(t)}=[a_{i,j}^{(t)}]_{1\ls i,j\ls 2n-1}$, we have
$$\per(A_n^{(t)})=\sum_{\sigma\in S_{2n-1}}\prod_{i=1}^{2n-1}(-1)^{[\![i<\sigma(i)\ \t{or}\ i=\sigma(i)<t]\!]}=\sum_{\sigma\in S_{2n-1}}(-1)^{\sum_{i=1}^{2n-1}[\![i<\sigma(i)\ \t{or}\ i=\sigma(i)<t]\!]},
$$
where we use the notation
$$[\![P]\!]=\begin{cases}1&\t{if the proposition}\ P\ \t{holds},
\\0&\t{otherwise}.\end{cases}$$

For each $t\in\{1,\ldots,2n\}$, we have
\begin{align*}\per(A_n^{(2n+1-t)})&=\sum_{\tau\in S_{2n-1}}(-1)^{\sum_{i=1}^{2n-1}[\![i<\tau(i)
\ \t{or}\ i=\tau(i)<2n+1-t]\!]}
\\&=\sum_{\tau\in S_{2n-1}}(-1)^{\sum_{i=1}^{2n-1}[\![2n-i>2n-\tau(i)
\ \t{or}\ 2n-i=2n-\tau(i)>t-1]\!]}
\\&=\sum_{\sigma\in S_{2n-1}}(-1)^{\sum_{j=1}^{2n-1}[\![j>\sigma(j)\ \t{or}\ j=\sigma(j)\gs t]\!]}
\\&=\sum_{\sigma\in S_{2n-1}}(-1)^{\sum_{j=1}^{2n-1}(1-[\![j\ls\sigma(j)\ \t{and}\ (j\not=\sigma(j)\ \t{or}\ j=\sigma(j)<t)]\!])}
\\&=-\sum_{\sigma\in S_{2n-1}}(-1)^{\sum_{j=1}^{2n-1}[\![j<\sigma(j)\ \t{or}\ j=\sigma(j)<t)]\!])}
=-\per(A_n^{(t)}).
\end{align*}
Therefore
$$\per(A_n)=\sum_{t=1}^{n}\l(\per(A_n^{(t)})+\per(A_n^{(2n+1-t)})\r)=0$$
as desired. \qed

\medskip
\noindent{\it Proof of Theorem \ref{Th1.2}}. Set $X_n=[x_{j,k}]_{1\ls j,k\ls 2n}$, where $x_{j,k}$
is defined as in \eqref{xjk}. Then $S(x_1,\ldots,x_{2n})=\per(X_n)$.

(i) We first prove that if one of the variables $x_1,\ldots,x_{2n}$ takes zero then $\per(X_n)=0$.
As $S(x_1,\ldots,x_{2n})$ is symmetric in $x_1,\ldots,x_{2n}$, without loss of generality we simply assume that $x_1=0$. We claim that
\begin{equation}\label{Xn}\per(X_n)=1+\sum_{k=1}^n\sum_{\tau\in C(2k)\atop 1\in D(\tau)}f(\tau)=0.
\end{equation}

Clearly, the permanent of
    $$
X_2=\begin{bmatrix}
	1 & -1\\
	1 & 1 \\
\end{bmatrix}
$$
vanishes,  and also
$$
1+f((1\ 2))=1+(-1)=0.
$$
Therefore \eqref{Xn} holds for $n=1$.

Now let $n>1$, and assume that \eqref{Xn} with $n$ replaced by any $m\in\{1,\ldots,n-1\}$
remains valid.

If $\tau$ is the identity $I$ of $S_{2n}$, then $f(\tau)=1$ since $D(\tau)=\emptyset$.
For any $\tau\in S_{2n}$ with $\tau\not=I$, we can decompose $\tau$ as a product of disjoint cycles:
$$(m_{11}\cdots m_{1l_1})(m_{21}\cdots m_{2l_2})\cdots (m_{k1}\cdots m_{kl_k}),$$
where $l_1,l_2,\ldots,l_k\in\{2,\ldots,2n\}$ and those $m_{il_i}\ (1\ls i\ls k,\ 1\ls j\ls l_i)$
are distinct numbers among $1,\ldots,2n$. For
$$\tau'=(m_{1l_1}\cdots m_{11})(m_{21}\cdots m_{2l_2})\cdots (m_{k1}\cdots m_{kl_k})$$
we have $f(\tau')=(-1)^{l_1}f(\tau)$ by Lemmas 2.1-2.2. Thus $f(\tau)+f(\tau')=0$ if $l_1$ is odd.
Therefore
\begin{equation}\label{perX}\per(X_n)=1+\sum_{\tau\in E_{2n}}f(\tau),
\end{equation}
where $E_{2n}$ consists of those $\tau\in S_{2n}$ which can be decomposed as a product
of distinct cycles of even length.

Let $\tau\in E_{2n}$ with $\tau(1)=1$. Then $|\Fix(\tau)|=2n-|D(\tau)|$ is even.
Write $|\Fix(\tau)|=2m$, where $1\ls m<n$. For any $\sigma\in C(2k)$ with $1\ls k\ls m$ and $D(\sigma)\cap D(\tau)=\emptyset$, we have $\tau\sigma\in E_{2n}$; if $1\ls j\ls 2n$ and $j\not\in \Fix(\tau)$ then $j\not\in D(\sigma)$ and hence $\sigma(j)=j$, thus we may view $\sigma$ as an element of $S_{2m}$. By the induction hypothesis, we have
$$1+\sum_{k=1}^m\sum_{\sigma\in C(2k),\,1\in D(\sigma)\atop D(\sigma)\cap D(\tau)=\emptyset}f(\sigma)=0.$$
Combining this with Lemma 2.1, we get
$$f(\tau)+\sum_{k=1}^m\sum_{\sigma\in C(2k),\,1\in D(\sigma)\atop D(\sigma)\cap D(\tau)=\emptyset}
f(\tau\sigma)=f(\tau)\(1+\sum_{k=1}^m\sum_{\sigma\in C(2k),\,1\in D(\sigma)\atop D(\sigma)\cap D(\tau)=\emptyset}f(\sigma)\)=0.$$

Any $\pi\in E_{2n}$ can be decomposed as a product of disjoint cycles of even length
 and thus we may write $\pi=\tau\sigma$, where $\sigma$ is the cycle containing $1$ in the decomposition of $\pi$ (if there are no cycle containing $1$ then $\sigma$ is the identity $I$ of $S_{2n}$), and $\tau$ is the product of all other cycles in the decomposition of $\pi$ (if there are no other cycles then $\tau$ is defined as the identity $I$ of $S_{2n}$).
 Note that $\tau(1)=1$, $1\in D(\sigma)$ and $D(\sigma)\cap D(\tau)=\emptyset$.
 Also, $\sigma\in C(2k)$ for some $1\ls k\ls m$.

 In view of \eqref{perX} and the last two paragraphs, we have
 \begin{equation}\label{per}\per(X_n)=1+\sum_{\pi\in E_{2n}}f(\pi)=1+\sum_{k=1}^{2n}\sum_{\sigma\in C(2k)\atop 1\in D(\sigma)}f(\sigma).
 \end{equation}
 This proves the first equality in \eqref{Xn}.

 Note that any $\tau\in C(2k)$ with $1\ls k\ls n$ and $1\in D(\tau)$ can be written as
 $(1\ m_2\ \cdots\ m_{2k})$, where $\{m_2,\ldots,m_{2k}\}$ is a subset of $\{2,3,\ldots,2n\}$
 with cardinality $2k-1$. Thus
 \begin{equation}\label{C(2k)}\sum_{\tau\in C(2k)\atop 1\in D(\tau)}f(\tau)=\bi{2n-1}{2k-1}s_k
 \end{equation}
 by Lemma 2.3, where  $s_1=s(0,x_2)=-1$, and
$s_k$ with $2\ls k\ls n$ denotes the constant $s(x_1,\ldots,x_{2k})$ in view of Lemma \ref{Lem2.3}.

 Combining \eqref{per} with \eqref{C(2k)}, we obtain
 \begin{equation} \label{per-C(2k)} \per(X_n)=1+\sum_{k=1}^n\bi{2n-1}{2k-1}s_k,\end{equation}
 which is independent of $x_2,\ldots,x_{2n}$. Therefore
 \begin{equation}\label{lim}\per(X_n)=\lim_{x_{2n}\to0}\cdots\lim_{x_2\to0}\per(X_n)
 =\sum_{\tau\in S_{2n}}\prod^{2n}_{j=1\atop\tau(j)\not=j}\lim_{x_n\to0}\cdots\lim_{x_2\to0}\f{x_j+x_{\tau(j)}}{x_j-x_{\tau(j)}}.
\end{equation}
If $1\ls j<k\ls 2n$, then
$$\lim_{x_{2n}\to0}\cdots\lim_{x_2\to0}\f{x_j+x_k}{x_j-x_k}
=\lim_{x_{2n}\to0}\cdots\lim_{x_{j+1}\to0}\f{0+x_k}{0-x_k}=-1.$$
If $1\ls k<j\ls 2n$, then
$$\lim_{x_{2n}\to0}\cdots\lim_{x_2\to0}\f{x_j+x_k}{x_j-x_k}
=\lim_{x_{2n}\to0}\cdots\lim_{x_{k+1}\to0}\f{x_j+0}{x_j-0}=1.$$
Thus, from \eqref{lim} we obtain
$$\per(X_n)=\sum_{\tau\in S_{2n}}\prod^{2n}_{j=1\atop\tau(j)\not=j}(-1)^{[\![j<\tau(j)]\!]}
=\sum_{\tau\in S_{2n}}\prod_{j=1}^{2n}a_{j,\tau(j)}=\per(A_n).$$
Combining this with Lemma 3.1, we finally get the equality $\per(X_n)=0$.
Therefore \eqref{Xn} holds.
This concludes our induction proof of the claim \eqref{Xn}.

(ii) By the proved claim \eqref{Xn} and \eqref{per-C(2k)}, for any integer $n>1$ we have
\begin{equation}\label{recurrence}
	1+\sum_{k=1}^{n}\binom{2n-1}{2k-1}s_k=0.
\end{equation}
Note that the last quality also holds for $n=1$ since $s_1=-1$.
Therefore, for $|x|<\pi/2$ we have
\begin{align*}
	0=&\sum_{n=1}^{\infty}\left(1+\sum_{k=1}^n\binom{2n-1}{2k-1}s_k\right)\f{x^{2n-1}}{(2n-1)!}\\
	=&\sinh x+\sum_{n=1}^{\infty}\(\sum_{k=1}^n\binom{2n-1}{2k-1}s_k\)\f{x^{2n-1}}{(2n-1)!}\\
	=&\sinh x+\sum_{k=1}^{\infty}s_k\f{x^{2k-1}}{(2k-1)!}\sum_{n=k}^{\infty}\f{x^{2n-2k}}{(2n-2k)!}\\
	=&\sinh x+\sum_{k=1}^{\infty}s_k\f{x^{2k-1}}{(2k-1)!}\cosh x
\end{align*}
and hence
$$
\sum_{k=1}^{\infty}(-1)^kT_k\f{x^{2k-1}}{(2k-1)!}=i\tan ix=-\tanh x=\sum_{k=1}^{\infty}s_k\f{x^{2k-1}}{(2k-1)!},
$$
which implies that $s_k=(-1)^kT(k)$ for all $k=1,2,3,\ldots$.
\medskip

In view of the above, we have completed the proof of Theorem \ref{Th1.2}.
 \qed

\begin{remark} For any integer $n\gs2$, as $s_1=-1$, from \eqref{recurrence} we obtain the equality
\begin{equation}\label{sk} \sum_{k=2}^n\bi{2n-1}{2k-1}s_k=2n-2.\end{equation}
\end{remark}

\section{Proof of Theorem 1.1}
\setcounter{lemma}{0}
\setcounter{theorem}{0}
\setcounter{equation}{0}
\setcounter{conjecture}{0}
\setcounter{remark}{0}
\setcounter{corollary}{0}

\noindent{\it Proof of Theorem \ref{Th1.1}}. We prove \eqref{S} by induction on $n$.

By \eqref{S(x12)}, we have 
\begin{equation}
\label{x1x2} S(x_1,x_2)=-\frac{4x_1x_2}{(x_1-x_2)^2}.
\end{equation}
Also, it is routine to verify that
$S(x_1,x_2,x_3,x_4)$ coincides with
$$16x_1x_2x_3x_4\left(\frac{1}{(x_1-x_2)^2(x_3-x_4)^2}
+\frac{1}{(x_1-x_3)^2(x_2-x_4)^2}+\frac{1}{(x_1-x_4)^2(x_2-x_3)^2}\right).$$
Thus \eqref{S} holds for $n=1,2$.

Now, let $n>2$ be an integer, and assume that \eqref{S} with $n$
replaced any $m\in\{1,\ldots,n-1\}$ remains valid.
By \eqref{perX},
\begin{equation}\label{S-sum}S(x_1,\ldots, x_{2n})=1+\sum_{\tau\in E_{2n}}f(\tau).\end{equation}
For any $\tau\in S_{2n}$, we define $d(\tau)$ as the least positive integer $d$ such that $\tau^d(1)=1$. Clearly $d(\tau)=1$ if and only if $1\in \Fix(\tau)$.
When $\tau\in E_{2n}$ and $1\not\in \Fix(\tau)$, obviously $d(\tau)$ is even. Thus
\begin{equation}\label{Sigma's}
	S(x_1,\ldots, x_{2n})=1+\Sigma_1+\Sigma_2+\Sigma_3,
\end{equation}
where
$$\Sigma_1=\sum_{\tau\in E_{2n}\atop d(\tau)=1}f(\tau),\ \
	\Sigma_2=\sum_{\tau\in E_{2n}\atop d(\tau)=2}f(\tau),\ \t{and}\
	\Sigma_3=\sum_{\tau\in E_{2n}\atop d(\tau)\gs 4}f(\tau).$$

For $2\ls i\ls 2n$, we write $x_2,\ldots,\hat x_i,\ldots,x_{2n}$
to denote the sequence $x_2,\ldots,x_{2n}$ with $x_i$ removed, thus
\begin{equation}\label{xi} S(x_2,\ldots, \hat x_i,\ldots,x_{2n})=1+\sum_{\tau\in E_{2n}
\atop 1,i\in\Fix(\tau)} f(\tau)\end{equation}
in view of \eqref{S-sum}.

We first deal with $\Sigma_1$.
By $\eqref{xi}$ we have
\begin{equation*}
	\sum_{i=2}^{2n}S(x_2,\ldots,\hat{x_i},\ldots, x_{2n})\\
   =2n-1+\sum_{\tau\in E_{2n}\atop 1\in \Fix(\tau)}f(\tau)\left(\#\Fix(\tau)-1\right),\\
\end{equation*}
and hence
\begin{equation}\label{Sigma-1}
	\Sigma_1=\sum_{i=2}^{2n}S(x_2,\ldots,\hat{x_i},\ldots, x_{2n})-2n+1-\sum_{\substack{\tau\in E_{2n}\\ 1\in \Fix(\tau)\\ \#\Fix(\tau)\gs4}}f(\tau)\left(\#\Fix(\tau)-2\right).
\end{equation}

Any $\tau\in E_{2n}$ with $d(\tau)=2$ can be decomposed as a product of
disjoint cycles of even length including the transpose $(1\ \tau(1))$.
Thus, with the aid of Lemma 2.1, we get
\begin{equation}\label{Sigma-2}
	\Sigma_2=\sum_{i=2}^{2n}f((1\ i))S(x_2,\ldots,\hat{x_i},\ldots, x_{2n}).
\end{equation}

Now we deal with $\Sigma_3$. Any $\tau\in E_{2n}$ with $d(\tau)\gs 4$
is a product of disjoint cycles, and hence we can write $\tau=\tau_1\tau_2$,
 where $\tau_1$ is a $d(\tau)$-cycle containing $1$, and $\tau_2$ is either the identity $I$
 or the product
 of all the other cycles in the decomposition of $\tau$.
 Thus, in view of Lemmas 2.1 and 2.3, we have
 $$\sum_{\tau\in E_{2n}\atop {d(\tau)\gs4\atop \tau_2=I}}f(\tau)
 =\sum_{k=2}^n\sum_{\tau_1\in C(2k)\atop 1\in D(\tau_1)}f(\tau_1)=\sum_{k=2}^n\bi{2n-1}{2k-1}s_k$$
 and
 \begin{align*}\sum_{\tau\in E_{2n}\atop {d(\tau)\gs4\atop \tau_2\not=I}}f(\tau)
 &=\sum_{\tau_2\in E_{2n}\atop {1\in\Fix(\tau_2)\atop\#\Fix(\tau_2)\gs4}}f(\tau_2)
 \sum_{k=2}^{\#\Fix(\tau_2)/2}\sum_{\tau_1\in C(2k)\atop 1\in D(\tau_1)}f(\tau_1)
 \\&=\sum_{\tau_2\in E_{2n}\atop {1\in\Fix(\tau_2)\atop\#\Fix(\tau_2)\gs4}}f(\tau_2)\sum_{k=2}^{\#\Fix(\tau_2)/2}
 \bi{\#\Fix(\tau_2)-1}{2k-1}s_k.
 \end{align*}
 Therefore
\begin{equation*}
	\Sigma_3=\sum_{k=2}^{n}\binom{2n-1}{2k-1}s_k+\sum_{\substack{\tau\in E_{2n}\\ 1\in \Fix(\tau)\\ \#\Fix(\tau)\gs4}}f(\tau)\sum_{k=2}^{\#\Fix(\tau)/2}\binom{\#\Fix(\tau)-1}{2k-1}s_k.
\end{equation*}
This, together with the identity \eqref{sk}, yields
\begin{equation}\label{Sigma-3}
	\Sigma_3=2n-2+\sum_{\substack{\tau\in E_{2n}\\ 1\in \Fix(\tau)\\ \#\Fix(\tau)\gs4}}f(\tau)\left(\#\Fix(\tau)-2\right).
\end{equation}
Combining \eqref{Sigma's} and \eqref{Sigma-1}-\eqref{Sigma-3}, we obtain
\begin{equation}\label{SS}
	S(x_1,\ldots, x_{2n})=\sum_{i=2}^{2n}(1+f((1\ i)))S(x_2,\ldots,\hat{x_i},\ldots, x_{2n}).
\end{equation}
For each $i=2,\ldots,2n$, clearly
$$
1+f((1\ i))=S(x_1,x_i)=-\frac{4x_1x_i}{(x_1-x_i)^2}
$$
by \eqref{x1x2},
and
$$S(x_2,\ldots,\hat{x_i},\ldots, x_{2n})=(-4)^{n-1}\prod_{j=2\atop j\not=i}^{2n}x_j
\times\sum_{\bigcup_{k=1}^{n-1}\{i_k,j_k\}=\{2,\ldots,2n\}\sm\{i\}}\prod_{k=1}^{n-1}\f1{(x_{i_k}-x_{j_k})^2}.$$
by the induction hypothesis. Thus, from \eqref{SS} we immediately get the desired \eqref{S}.

 In view of the above, we have completed our proof of \eqref{S} by induction. \qed

\section{Proof of Theorem 1.3}
\setcounter{lemma}{0}
\setcounter{theorem}{0}
\setcounter{equation}{0}
\setcounter{conjecture}{0}
\setcounter{remark}{0}
\setcounter{corollary}{0}

\begin{lemma}[Guo et al. \cite{G}]\label{known-lemma}
	For any integer $n>2$ and variables $x_1,\ldots,x_n$, we have
	$$
	\sum_{\tau\in C(n)}\prod_{j=1}^{n}\f{1}{x_{\tau(j)}-x_j}=0.
	$$
\end{lemma}

For any positive integer $n$, let ${\mathcal D}(n)$ denote the set of all derangements of $1,\ldots,n$.
The following result was originally conjectured by Sun \cite{ConjBook,S1} and recently proved by X. Guo, X. Li, Z. Tao and T. Wei \cite{G}.

\begin{lemma}
	Let $n>1$ be an integer,  and let $\zeta$ be a primitive $n$th root of unity.
	
{\rm (i)} If $n$ is even, then
		\begin{equation}\label{known-identity1}
			\sum_{\tau\in {\mathcal D}(n)}\prod_{j=1}^{n}\f{1}{1-\zeta^{j-\tau(j)}}=\f{((n-1)!!)^2}{2^n}.
		\end{equation}

{\rm (ii)} If $n$ is odd, then
	    \begin{equation} \label{known-identity2}
	    	\sum_{\tau\in {\mathcal D}(n-1)}\prod_{j=1}^{n}\f{1}{1-\zeta^{j-\tau(j)}}=\f{1}{n}\left(\f{n-1}{2}!\right)^2.
	    \end{equation}
\end{lemma}

\noindent{\it Proof of Theorem \ref{Th1.3}}. (i) For $1\ls j,k\ls 2n$, clearly
$$\f{x_{j,k}-1}{2x_k}=\begin{cases}\f{1}{x_j-x_k}&\text{if}\ j\not=k,\\0&\text{if}\ j=k.\end{cases}$$
By Lemma \ref{known-lemma}, we have
\begin{align*}
	\per\left[\f{x_{j,k}-1}{2x_k}\right]_{1\ls j,k\ls 2n}&=\sum_{\tau\in {\mathcal D}(2n)}\prod_{j=1}^{2n}\f{1}{x_j-x_{\tau(j)}}
=\sum_{\tau\in {\mathcal D}(2n)\atop \tau^2=I}\prod_{j=1}^{2n}\f{1}{x_j-x_{\tau(j)}}
\\&=\sum_{\bigcup_{k=1}^n\{i_k,j_k\}=\{1,\dots,2n\}} \prod_{k=1}^n\frac{1}{(x_{i_k}-x_{j_k})(x_{j_k}-x_{i_k})},
\end{align*}
and
\begin{align*}
	\det\left[\f{x_{j,k}-1}{2x_k}\right]_{1\ls j,k\ls 2n}&=\sum_{\tau\in {\mathcal D}(2n)}\sgn(\tau)\prod_{j=1}^{2n}\f{1}{x_j-x_{\tau(j)}}
=\sum_{\tau\in {\mathcal D}(2n)\atop \tau^2=I}\sgn(\tau)\prod_{j=1}^{2n}\f{1}{x_j-x_{\tau(j)}}
\\
	&=\sum_{\bigcup_{k=1}^n\{i_k,j_k\}=\{1,\dots,2n\}} (-1)^n\prod_{k=1}^n\frac{1}{(x_{i_k}-x_{j_k})(x_{j_k}-x_{i_k})}.
\end{align*}
Combining these with \eqref{S}, we immediately obtain the desired \eqref{per-det identity}.

(ii) Now we turn to prove part (ii) of Theorem \ref{Th1.3}. For $1\ls j,k\ls n$, clearly
$$
c_{j,k}-1=\begin{cases}\f{2\zeta^{j-k}}{1-\zeta^{j-k}} &\t{if $j\not=k$,}\\
	0&\t{if $j=k$.}\end{cases}
$$

{\it Case} 1. $n$ is even.

Let $x_j=\zeta^j$ for $j=1,\ldots,n$. Then $c_{j,k}=x_{k,j}$ for all $j,k=1,\ldots,n$.
 Thus, by \eqref{per-det identity} and \eqref{known-identity1}, we have
	\begin{align*}
	    \per[c_{j,k}]_{1\ls j,k\ls n}&=\per[c_{j,k}-1]_{1\ls j,k\ls n}
	    =2^n\prod_{j=1}^n\zeta^j\times\prod_{k=1}^n\zeta^{-k}\times\sum_{\tau\in {\mathcal D}(n)}\prod_{j=1}^{n}\f{1}{1-\zeta^{j-\tau(j)}}\\
	    &=2^n\times\f{((n-1)!!)^2}{2^n}=((n-1)!!)^2.
	\end{align*}
This proves \eqref{c-even}.

    {\it Case} 2. $n$ is odd.

    In light of \eqref{per-det identity} and \eqref{known-identity2}, we have
 	\begin{align*}
 	\per[c_{j,k}]_{1\ls j,k\ls n-1}&=\per[c_{j,k}-1]_{1\ls j,k\ls n-1}
 	=2^{n-1}\sum_{\tau\in {\mathcal D}(n-1)}\prod_{j=1}^{n}\f{1}{1-\zeta^{j-\tau(j)}}
 	\\&=2^{n-1}\times\f{1}{n}\left(\f{n-1}{2}!\right)^2
 	=\f{((n-1)!!)^2}{n}.
    \end{align*}
    This proves \eqref{c-odd}.

     In view of the above, we have completed our proof of Theorem 1.3. \qed


\begin{thebibliography}{99}

\bibitem {FLS} S. Fu, Z. Lin and Z.-W. Sun, {\it Proof of several conjectures relating permanents to combinatorial sequences},  arXiv:2109.11506, preprint, 2021.

\bibitem{G} X. Guo, X. Li, Z. Tao, T. Wei, {\it The eigenvectors-eigenvalues identity and Sun's conjectures on determinants and permanents} , arXiv:2206.02592, preprint, 2022.

\bibitem{Sloane} N. J. A. Sloane, Sequence A000182 at OEIS (The On-Line Encyclopedia of Integer Sequences), {\tt http://oeis.org/A000182}.


\bibitem{ConjBook} Z.-W. Sun, New Conjectures in Number Theory and Combinatorics (in Chinese),
Harbin Institute of Technology Press, Harbin, 2021.

\bibitem{S1} Z.-W. Sun, {\it Arithmetic properties of some permanents}, arXiv:2108.07723, preprint, 2022.

\bibitem{S2} Z.-W. Sun, {\it On some determinants and permanents}, arXiv:2207.13039, preprint, 2022.

 \bibitem{OEIS} Z.-W. Sun, Sequence A356041 at OEIS (The On-Line Encyclopedia of Integer Sequences),
2002. {\tt http://oeis.org/A356041}

\bibitem{Sun} Z.-W. Sun, {\it A conjectural permanent identity}, Question 427232 at MathOverflow,
July 24, 2022. {\tt https://mathoverlow.net/questions/427232}

\bibitem{WangSun} H. Wang and Z.-W. Sun, {\it Proof of a conjecture involving derangements and roots of unity}, arXiv:2206.02589, preprint, 2022.



\end{thebibliography}
\end{document}